\documentclass[aos,seceqn,nameyear,dvips]{arximspdf}
\usepackage{dcolumn}

\doi{10.1214/09-AOS758}
\volume{38}
\issue{3}
\pubyear{2010}
\firstpage{1546}
\lastpage{1567}

\makeatletter

\newcolumntype{d}[1]{D{.}{.}{#1}}

\newtheorem{theorem}{Theorem}[section]
\newtheorem{crl}{Corollary}[section]
\newtheorem{lem}{Lemma}[section]
\newtheorem{prp}{Proposition}[section]

\newproclaim{Remark}{Remark}

\makeatother

\begin{document}
\begin{frontmatter}

\title{Exact properties of Efron's biased coin randomization procedure}
\runtitle{Biased coin randomization}

\begin{aug}
\author[A]{\fnms{Tigran} \snm{Markaryan}\ead[label=e1]{tigranmarkaryan@yahoo.com}} and
\author[B]{\fnms{William F.} \snm{Rosenberger}\corref{}\thanksref{t1}\ead[label=e2]{wrosenbe@gmu.edu}}
\runauthor{T. Markaryan and W. F. Rosenberger}
\affiliation{George Mason University}
\address[A]{2134 Glencourse Lane\\
Reston, Virginia 20191\\
USA\\
\printead{e1}}
\address[B]{Department of Statistics\\
George Mason University\\
4400 University Drive, MS 4A7\\
Fairfax, Virginia 22030-4444\\
USA\\
\printead{e2}}
\end{aug}

\thankstext{t1}{Supported by NSF Grant DMS-09-04253 under the
2009 Recovery and Reinvestment Act.}

\received{\smonth{6} \syear{2009}}
\revised{\smonth{9} \syear{2009}}

%
\begin{abstract}
Efron [\textit{Biometrika} \textbf{58} (1971) 403--417] developed a
restricted randomization procedure to promote balance between two
treatment groups in a sequential clinical trial. He called this the
\textit{biased coin design}. He also introduced the concept of
\textit{accidental bias}, and investigated properties of the procedure
with respect to both accidental and selection bias, balance, and
randomization-based inference using the steady-state properties of the
induced Markov chain.
In this paper we revisit this procedure, and derive closed-form
expressions for the exact properties of the measures derived
asymptotically in Efron's paper. In particular, we derive the exact
distribution of the treatment imbalance and the variance-covariance
matrix of the treatment assignments. These results have application in
the design and analysis of clinical trials, by providing exact formulas
to determine the role of the coin's bias probability in the context of
selection and accidental bias, balancing properties and
randomization-based inference.
\end{abstract}

%

\end{frontmatter}

\section{Introduction}

Efron (\citeyear{Efron1971}) introduced his famous \textit{biased coin
design} as a method that ``\ldots tends to balance the experiment, but
at the same time is not over vulnerable to various common forms of
experimental bias.'' The primary application is in sequential clinical
trials where balance in the numbers randomly assigned to two treatment
groups is sometimes desirable for power considerations. In such cases,
it is also desirable to maintain near-balance at intermediate points in
the trial as heterogeneity or time trends in patient characteristics
may lead to less comparable treatment arms. Randomization protects from
imbalances in unknown covariates related to outcomes (which Efron
referred to as \textit{accidental bias}, introduced for the first time
in the 1971 paper), selection bias and provides a basis for inference.
Efron explored the balancing properties of the biased coin design, as
well as its susceptibility to selection and accidental bias, and
discussed the implications for randomization-based inference. All of
these results were based on studying the steady-state properties of the
Markov chain induced by the imbalance process of biased coin
randomization.

Let $\mathbf{T}_n=(T_1,\ldots,T_n)'$ be a randomization sequence, where
$T_i=1$ if treatment $A$ is assigned, and $T_i=-1$ if treatment $B$ is
assigned, $i=1,\ldots,n$. After $j$ assignments, let $D_{j}$ be the
difference in the number of patients assigned to treatments $A$ and
$B$; that is, $D_j=\sum_{i=1}^j T_i$. The biased coin design with bias
$p\in [0.5,1]$, denoted $\mathit{BCD}(p)$, is defined by
\begin{eqnarray*}
P(T_{j}=1) = \cases{
1/2, &\quad when $D_{j-1} = 0$, \cr
p, &\quad when $D_{j-1} < 0$,\qquad$j = 1,2,3,\ldots,$\cr
1-p, &\quad when $D_{j-1} > 0$.}
\end{eqnarray*}
Note that $p=0.5$ results in complete randomization and $p=1$ results
in a permuted block design with block size of 2, in which case every
alternate assignment is deterministic. Efron notes that the
$\{|D_{n}|\}_{n=1}^{\infty}$ process forms a Markov chain of period 2
with states $0, 1, 2,\ldots$ and a reflecting barrier at the origin. He
then proves that the $|D_{n}|$ process has stationary probabilities
$\pi_{j}$, given by
%
%
\begin{equation}
\label{eq:def_BCD_stationaryDist}
\pi_{j} = \cases{
\dfrac{r^2-1}{2r^{j+1}}, &\quad when $j \geq1$, \vspace*{2pt}\cr
\dfrac{r-1}{2r}, &\quad when $j = 0$,}
\end{equation}
where $q=1-p$, $r=p/q \geq1$. Efron uses the formulas obtained for
stationary probabilities to write the form of the limiting
probabilities of perfect balance ($n$ is even) and imbalance of 1 ($n$
is odd) as
\begin{eqnarray*}
\lim_{n \rightarrow\infty} P(|D_{2m}| = 0) & = & 2\pi
_{0} = \frac{r-1}{r}, \\
\lim_{n \rightarrow\infty} P(|D_{2m+1}| = 1) & = & 2\pi
_{1} = \frac{r^2-1}{r^2}.
\end{eqnarray*}

Most research on the theory of randomization in recent years has
focused on generalizations of Efron's procedure [see, e.g., Wei
(\citeyear{Wei1978}),
Soares and Wu (\citeyear{SoaresWu1982}), Eisele (\citeyear{Eisele94}),
Chen (\citeyear{Chen1999}), Baldi Antognini and Giovagnoli (\citeyear
{BaldiAntogniniGiovagnoli2004}) and Hu and Zhang
(\citeyear{HZ04})] rather than Efron's procedure itself. In
particular, Baldi
Antognini and Giovagnoli's (\citeyear{BaldiAntogniniGiovagnoli2004})
``adjustable biased
coin design'' is stochastically more balanced, and therefore uniformly more
powerful, than the other procedures [Baldi Antognini (\citeyear
{BaldiAntognini08})].

The remainder of Efron's article is devoted to selection bias, as
defined by Blackwell and Hodges (\citeyear{BlackwellHodges1957}),
accidental bias and
randomization as a basis for inference. Efron notes that the best
guessing strategy against the $\mathit{BCD}(p)$ is to always guess the group
that has occurred least often up to that point. The probability of
correctly guessing at the $j$th step is
%
%
\begin{equation}
\label{eq:BCD_SelectionBias}
\tfrac{1}{2}P(D_{j-1} = 0)+ p P(|D_{j-1}| > 0),
\end{equation}
which asymptotically approaches $1/2 + (r-1) / 4r$ and therefore has
asymptotic excess selection bias of
%
%
\begin{equation}\label{selbias}
\frac{(r-1)}{4r}.
\end{equation}

Accidental bias refers to the squared bias of the treatment effect in a
linear regression when an unknown covariate $\mathbf{z}$ is
left out of the model. Efron derives this bias as
\[
E(\mathbf{z}' \mathbf{T}_{n})^2 = \mathbf{z}' \bolds
{\Sigma}_{\mathbf{T}_n} \mathbf{z},
\]
where $\bolds{\Sigma}_{\mathbf{T}_n}=\operatorname{Var}(\mathbf{T}_n)$.
He suggests a minimax approach by noting that
%
%
\begin{equation}
\label{minimax}
\mathbf{z}' \bolds{\Sigma}_{\mathbf{T}_n} \mathbf{z}
\leq
\mbox{maximum eigenvalue of } \bolds{\Sigma}_{\mathbf{T}_n},
\end{equation}
where the inequality follows from the assumption that $\|\mathbf{z}\|=1$.
Note that the minimum possible value for the maximum
eigenvalue is 1 which corresponds to complete randomization.
Instead of directly examining $\bolds{\Sigma}_{\mathbf{T}_n}$
(which he acknowledges is difficult), Efron looks at the much simpler
process $T_1, T_2, T_3,\ldots,T_n$, assuming that it is stationary, and
aims at finding the asymptotic covariance structure of the process. He
then shows that the asymptotic maximum eigenvalue of the covariance
vector $(T_{h+1}, \ldots, T_{h+N})$ as $h \rightarrow\infty$,
$\lambda
_N$, is increasing in $N$ and has a finite limit. Based on numerical
evidence, Efron conjectures that $\lim_{N \rightarrow\infty} \lambda_N
= 1+(p-q)^2$. This was later proved by Steele (\citeyear{Steele1980}).
However, Smith
(\citeyear{Smith1984}) shows by counterexample that Efron's solution
may be
unsatisfactory when there are short-term dependencies in the data.

In this paper, we derive exact properties of Efron's procedure. In
particular, in Section \ref{sec2}, we derive a closed-form expression for the
distribution of $D_n$ and give the
explicit form of $\bolds{\Sigma}_{\mathbf{T}_n}$. These
formulas are remarkably compact for the complexity of the problems. We
describe computational considerations in Section \ref{sec3}.
In Section \ref{sec4}, we apply these results to deriving an explicit form for
the excess selection bias, prove a result on the maximum eigenvalue of
$\bolds{\Sigma}_{\mathbf{T}_n}$
and discuss randomization as a basis for inference. We also compare the
exact results with Efron's for various $n$ and $p$. In Section \ref{sec5}, we
draw conclusions. Finally, all proofs are given
in Appendices \ref{appA}--\ref{appC}.

\section{Exact distribution of $D_n$ and ${\Sigma}_{\mathbf{T}_n}$}\label{sec2}

We will assume the following conventions throughout the mathematical
developments.
\begin{enumerate}
\item For brevity, we adopt the convention to treat a combination
${x \choose y}$ as zero whenever any of the following conditions is
true: $x <
0, y <0, x < y, y$ is not an integer.
\item We treat summations as 0 if the upper limit of the summation is
smaller than the lower limit.
\item We treat conditional probabilities, conditional on
zero-probability events to be~0.
\end{enumerate}

The distribution of $D_n$ requires determination of the exact
distribution of a denumerable homogeneous random walk.
The following result is given as the first theorem:
\begin{theorem}
\label{thm:distdnshort}
Let $n = 1,2,3,\ldots, 0 \leq k \leq n$ and $n$ and $k$ have
the same parity.
Then, the distribution of $D_{n}$ of the $\mathit{BCD}(p)$ is given by formulas
(\ref{eq:DistDnShortKpos}) and (\ref{eq:DistDnShortKzero}).

For $k>0$,
%
%
\begin{equation}
\label{eq:DistDnShortKpos}\quad
P(D_{n}=\pm k)= \frac{1}{2} p^{({n-k})/{2}}
\sum_{l=0}^{({n-k})/{2}}
{\frac{n+k-2l}{n+k+2l}
\pmatrix{\dfrac{n+k}{2}+l \vspace*{2pt}\cr
l}q^{k+l-1}.}
\end{equation}
For $k=0$,
%
%
\begin{equation}
\label{eq:DistDnShortKzero}
P(D_{n}=0)=p^{{n}/{2}}
\sum_{l=0}^{{n}/{2}-1}
{\frac{n-2l}{n+2l}
\pmatrix{
\dfrac{n}{2}+l \vspace*{2pt}\cr
l}q^{l}.}
\end{equation}
\end{theorem}
\begin{pf}
See Appendix \ref{appA}.
\end{pf}

The compact form of these equations arises from patterns in
polynomials of $p$ and $q$ that can be seen
developing for small $n$ as $n$ increments. The proof is then by
induction. Note that the distribution of $N_A(n)$ follows immediately,
since $N_A(n)=(D_n+n)/2$.

Define $t_k=P(T_n=1|D_{n-1}=k)$. We now derive the
covariance of $(T_n,T_m)$.
\begin{theorem}\label{thm:distdndm}
Let $1 \leq n < m$. Then the joint distribution of $(T_{n},
T_{m})$ of the $\mathit{BCD}(p)$, $p \in[1/2,1]$, is given by
%
%
\begin{equation}
\label{eq:JointTnTm}\qquad
P(T_{n}=1, T_{m}=1) =
\sum_{k=-n+1}^{n-1}
{\biggl(
\biggl(
\frac{1}{2}-t_{k+1}
\biggr)
\hat{f}_{k+1,0}^{(m-n-1)} + t_{k+1}
\biggr)
d_{n-1,k} t_{k}},
\end{equation}
where
%
%
\[
d_{n,k} = P(D_{n}=k) \mbox{ and is given in (\ref{eq:DistDnShortKpos})
and (\ref{eq:DistDnShortKzero})}
\]
and
\begin{equation}
\label{def:f_hat}\qquad
\hat{f}_{k,0}^{(u)} = \sum_{l=|k|}^{u} {f_{k,0}^{(l)}}=
\cases{\displaystyle\sum_{l=|k|}^{u} {
\frac{|k|}{l}
\pmatrix{l \vspace*{2pt}\cr\dfrac{l+|k|}{2}}
p^{({l+|k|})/{2}} q^{({l-|k|})/{2}}
},&  when $k \neq0$;\hspace*{-6pt}\cr
1, & when $k = 0$.\hspace*{-6pt}}
\end{equation}
\end{theorem}
\begin{pf}
See Appendix \ref{appB}.
\end{pf}

The form of $\bolds{\Sigma}_{\mathbf{T}_n}$ follows
immediately:
\begin{crl}
\label{crl:SigmaDn}
Let $\bolds{\Sigma}_{\mathbf{T}_n}$ be the
covariance matrix of $\mathbf{T}_n$ of the $\mathit{BCD}(p)$, $p \in[1/2,1]$.
Then the $(i,j)$th entry of the matrix, $\sigma_{ij}$, where $1 \leq i
\leq j \leq n$, is given by
%
%
\begin{equation}
\label{eq:Sigma_i_j}
\sigma_{ij}=
\cases{\displaystyle 4 \cdot
\sum_{k=-i+1}^{i-1}
 \biggl(
\biggl(
\frac{1}{2}-t_{k+1}
\biggr)
\hat{f}_{k+1,0}^{(j-i-1)} + t_{k+1}
\biggr)\cr
\hspace*{45pt}{}\times
d_{i-1,k} t_{k} - 1
,\qquad\mbox{when $i < j$}; \cr
1, \hspace*{127pt} \mbox{when $i = j$};}
\end{equation}
where $\hat{f}_{k,0}^{(u)}$ is defined in (\ref{def:f_hat}).
\end{crl}

\section{Computational considerations}\label{sec3}

This section contains some observations on the computation of $P(D_n =
k)$ according to formulas (\ref{eq:DistDnShortKpos}) and (\ref
{eq:DistDnShortKzero}) and the computation of $P(T_{n}=1, T_{m}=1)$
according to formulas (\ref{eq:JointTnTm}) and (\ref{def:f_hat}). These
formulas involve terms that are products of large factorials and powers
of numbers that are between $0$ and $1$. The key is to calculate these
products in such order that the result does not get too large or too
small too quickly. We focus on the computation of (\ref
{eq:DistDnShortKpos}) here as the other formulas are similar. For $n
\le100$, calculating the combination and multiplying by powers of $p$
and $q$ directly works well. However, for larger values of $n$,
precision may be lost if the intermediate products become too large or
too small.

Formula (\ref{eq:DistDnShortKpos}) involves $(n-k)/2 + 1$ terms, each
of which is a product of powers of $p$, powers of $q$, positive
integers and reciprocals of positive integers. The generic term of the
right-hand side of (\ref{eq:DistDnShortKpos}) can be written as
%
%
\begin{equation}
\label{eq:Calc_of_Dnk}\qquad
\frac{1}{2} \frac{n+k-2l}{n + k + 2l}
\overbrace{p \cdots p}^{({n - k})/{2}}
\cdot\overbrace{q \cdots q}^{k + l -1}
\cdot\overbrace{\frac{1}{2} \cdot\frac{1}{3} \cdots\frac
{1}{l}}^{l - 1}
\cdot\overbrace{ \biggl(\frac{n + k}{2} + 1 \biggr) \cdots
\biggl(\frac{n + k}{2} + l \biggr)}^{l}.
\end{equation}
There are $(n + k) / 2 + 2l$ factors less than 1 and $l$
greater than 1. Denote these two groups by $\{a_s\}_{s=1}^{(n + k)/2 +
2l}$ and $\{b_s\}_{s=1}^{l}$, respectively. Assume that the $a_s$ are
indexed in decreasing order.
The following simple algorithm ensures that the running products for
calculating (\ref{eq:Calc_of_Dnk}) do not become too large or too small
early on in the calculation process.
\begin{enumerate}
\item Fix a number $M$ that the running product cannot exceed. Any
number that is larger than $2n$ will work.
\item Fix a number $m$ that is close to the machine epsilon. When the
running product gets close to $m$, the algorithm will know that further
multiplication by small numbers may result in loss of precision.
\item \hypertarget{mult_large} Start multiplication using numbers in $\{b_s\}
_{s=1}^{l}$ until the running product exceeds~$M$.
\item \hypertarget{mult_small} Multiply the running product with numbers in
$\{a_s\}_{s=1}^{(n + k) / 2 + 2l}$ until the running product is less
than~$M$.
\item Iterate through Steps \hyperlink{mult_large}{3} and \hyperlink{mult_small}{4}
until numbers in $\{b_s\}_{s=1}^{l}$ are depleted.
\item Continue multiplying the running product with the remaining
numbers in $\{a_s\}_{s=1}^{(n + k)/2 + 2l}$, from the largest to the
smallest. Two cases are possible:
the final product is larger than $m$, in which case the algorithm is completed;
at some point, the running product becomes less than $m$, in which
case one can save the result as a product of two or more small numbers.
\end{enumerate}
For example, we used the algorithm to compute Table \ref
{tab:Dnk_errors_within_smallk}, which gives the value of $n$ at which
the steady state probabilities are within certain
percentages of the exact probability $P(D_n=k)$, for various $k$ and
$p$. The same idea can be used for calculating (\ref{def:f_hat}).

%
\begin{table}
\tablewidth=250pt
\caption{Values of $n$ starting at which, steady state probabilities
are within~10\%, 5\%, 1\% and 0.1\% of $P(D_n = k)$, $k = 0, 1, 2, 3, 4$}
\label{tab:Dnk_errors_within_smallk}
\begin{tabular*}{\tablewidth}{@{\extracolsep{\fill}}lcd{3.0}d{3.0}d{3.0}d{4.0}@{}}
\hline
& & \multicolumn{4}{c@{}}{\textbf{Errors within}}\\[-4pt]
& & \multicolumn{4}{c@{}}{\hrulefill}\\
$\bolds k$ & $\bolds p$ & \multicolumn{1}{c}{\textbf{10\%}}
& \multicolumn{1}{c}{\textbf{5\%}} & \multicolumn{1}{c}{\textbf{1\%}}
& \multicolumn{1}{c@{}}{\textbf{0.1\%}}\\
\hline
\phantom{0}0 & 0.6 & 20 & 34 & 74 & 146 \\
& 0.7 & 6 & 8 & 18 & 34 \\
& 0.8 & 2 & 4 & 8 & 14 \\
& 0.9 & 2 & 2 & 4 & 6 \\
[3pt]
\phantom{0}1 & 0.6 & 19 & 33 & 73 & 145 \\
& 0.7 & 5 & 7 & 17 & 33 \\
& 0.8 & 1 & 3 & 7 & 13 \\
& 0.9 & 1 & 1 & 3 & 5 \\
[3pt]
\phantom{0}2 & 0.6 & 14 & 28 & 68 & 140 \\
& 0.7 & 4 & 4 & 8 & 22 \\
& 0.8 & 4 & 4 & 8 & 14 \\
& 0.9 & 2 & 4 & 6 & 8 \\
[3pt]
25 & 0.6 & 183 & 211 & 279 & 379 \\
& 0.7 & 85 & 93 & 113 & 141 \\
& 0.8 & 53 & 57 & 65 & 77 \\
& 0.9 & 37 & 39 & 43 & 49 \\
[3pt]
50 & 0.6 & 342 & 380 & 464 & \mbox{$>$}500 \\
& 0.7 & 158 & 168 & 194 & 226 \\
& 0.8 & 100 & 104 & 116 & 130 \\
& 0.9 & 70 & 72 & 78 & 86 \\
\hline
\end{tabular*}
\end{table}

Finally, for the computation of $\bolds{\Sigma}_{\mathbf{T}_n}$, the
following proposition gives a property of the matrix
that can facilitate computation. The proof follows from Corollary \ref{crl:SigmaDn}
and Lemma \ref{lem:Joint_To_Conditional} and is omitted.
\begin{prp}
If $\bolds{\Sigma}_{\mathbf{T}_n}$ is partitioned into $2
\times2$ submatrices, then all the off-diagonal
submatrices are constant (i.e., have the same elements in both rows and
columns).
\end{prp}

\section{Applications to clinical trials}\label{sec4}

In this section we apply the results of Section \ref{sec2} to the study of
balancing properties of the randomization procedure, selection and
accidental biases and randomization
as a basis for inference. Each of these is a consideration in the
appropriate selection of a randomization procedure in clinical trials
[see Rosenberger and Lachin (\citeyear{RosenbergerLachin2002})].

\subsection{Balancing properties of the biased coin design}
All finite balancing properties of the biased coin design can be
investigated with the help of Theorem~\ref{thm:distdnshort} which
provides the means for exact calculations of the probabilities
involving $P(D_n = k)$. In particular, the exact variance is given in
the following proposition:
\begin{prp}
The exact variance of $D_n$ is given by
%
%
\begin{equation}
\label{eq:VarOfDn}\qquad
\operatorname{Var}(D_{n}) =
\mathop{\sum_{k = 1}}_{n-k\ \mathrm{even}}^{n}{k^2
p^{({n-k})/{2}}
\sum_{l=0}^{({n-k})/{2}}
{\frac{n+k-2l}{n+k+2l}
\pmatrix{
\dfrac{n+k}{2}+l \vspace*{2pt}\cr
l}q^{k+l-1}.}}
\end{equation}
\end{prp}

The variance of the imbalance of the biased coin design for
different values of $n$ and $p$ is provided in Table \ref
{tab:VarianceOfImbalanceBCD}.
%
%
\begin{table}[b]
\tablewidth=300pt
\caption{Variance of the imbalance of the $\mathit{BCD}(p)$ for different
values of $n$ and $p$}\label{tab:VarianceOfImbalanceBCD}
\begin{tabular*}{\tablewidth}{@{\extracolsep{\fill}}ld{2.2}ccc@{}}
\hline
& \multicolumn{4}{c@{}}{$\bolds p$} \\[-4pt]
& \multicolumn{4}{c@{}}{\hrulefill}\\
 & \multicolumn{1}{c}{\textbf{0.6}}
& \multicolumn{1}{c}{\textbf{0.7}} & \multicolumn{1}{c}{\textbf{0.8}}
& \multicolumn{1}{c@{}}{\textbf{0.9}} \\
\hline
\multicolumn{5}{c}{$n$ even}\\[4pt]
\phantom{0}10 & 5.19 & 2.55 & 1.18 & 0.46 \\
\phantom{0}20 & 7.65 & 2.91 & 1.21 & 0.46 \\
\phantom{0}50 & 10.78 & 3.04 & 1.21 & 0.46 \\
100 & 12.10 & 3.04 & 1.21 & 0.46 \\
200 & 12.45 & 3.04 & 1.21 & 0.46 \\
$\infty$ & 12.48 & 3.04 & 1.21 & 0.46 \\
[4pt]
\multicolumn{5}{c}{$n$ odd}
\\
[4pt]
\phantom{0}5 & 3.30 & 2.15 & 1.45 & 1.10 \\
15 & 6.63 & 2.95 & 1.56 & 1.10 \\
25 & 8.52 & 3.13 & 1.57 & 1.10 \\
75 & 11.73 & 3.20 & 1.57 & 1.10 \\
$\infty$ & 12.52 & 3.21 & 1.57 & 1.11 \\
\hline
\end{tabular*}
\end{table}
Also given in the table is the limiting variance based on the steady
state distribution of the induced Markov chain.
The formulas for odd and even $n$ are given in the following
proposition which follows directly from (\ref{eq:def_BCD_stationaryDist}).
\begin{prp}
\label{prp:AsymptoticVarianceDn}
Under the limiting distribution of the $\mathit{BCD}(p)$, $p \in
[1/2,1]$, the variance of the imbalance is given by
%
%
\begin{eqnarray}
&\displaystyle\frac{4r(r^2+1)}{(r^2-1)^2}\qquad\mbox{when number of trials is
even},&
\nonumber\\[-8pt]\\[-8pt]
&\displaystyle\frac{8r^2}{(r^2-1)^2} + 1\qquad\mbox{when number of trials is
odd.}&\nonumber
\end{eqnarray}
\end{prp}

As can be seen in the table, odd and even $n$ form different patterns.
This is due to the differences in the supports of the distributions; in
particular, a significant mass is concentrated at $0$ when $n$ is even
and $p$ is large. Note that both odd and even $n$ form an increasing
series for each $p$. This is expected and follows from Theorem 1 in
Efron (\citeyear{Efron1971}) with $h(j) = j^2$. It is also the case
that $\operatorname{Var}(D_{n})$
is a decreasing function of $p$ for each $n$. This is also expected and
was proved by Efron as Theorem 3 with $h(j) = j^2$. It is clear that
balancing properties stabilize for moderate-sized trials of around 75
to 100. This contrasts to
other randomization procedures such as complete randomization and the
urn design [Wei (\citeyear{Wei1978})] where the variance of $D_n$ grows
at a rate
$O(n)$ [Rosenberger and Lachin (\citeyear{RosenbergerLachin2002}),
Chapter~3)].

\subsection{Selection bias}
Theorem \ref{thm:distdnshort} allows us to calculate the selection bias
for the $\mathit{BCD}(p)$ using (\ref{eq:BCD_SelectionBias}).
When $n$ is even, $P(D_{n-1} = 0) = 0$, and therefore the selection
bias is $p$.
Obviously, the selection bias when $n = 1$ is $1/2$.
When $n$ is an odd number exceeding $1$, $n = 2m + 1$ and $m \in
\mathbb
{N}$, substituting the expression for $P(D_{n-1} = 0)$
from Theorem \ref{thm:distdnshort}, we obtain the following expression
for the selection bias for this case:
%
%
%
\begin{equation}
\label{eq:SelectionBiasBCD}
p - \biggl(p - \frac{1}{2}\biggr)
p^{m}
\sum_{l=0}^{m - 1}
{
\frac{m - l}{m + l}
\pmatrix{m + l \cr
l}q^{l}.
}
\end{equation}
Now we can formulate a result on the total selection bias in
$n$ trials.
\begin{prp}
\label{thm:SelectionBiasBCD}
The total amount of selection bias in $n$, $n \ge1$, trials
for the $\mathit{BCD}(p)$ is given by
%
%
\begin{equation}
\label{eq:TotalSelectionBiasBCD}
\frac{1}{2} + (n - 1)p -
\biggl(p - \frac{1}{2}\biggr)
\sum_{m = 1}^{[(n - 1)/2]}
p^{m}
\sum_{l=0}^{m - 1}
{
\frac{m - l}{m + l}
\pmatrix{m + l \cr
l}q^{l}},
\end{equation}
where $[a]$ denotes the integer part of $a$ and we use the
adopted convention that the sum is treated as zero when the upper limit
of summation is smaller than the lower limit.
\end{prp}

One subtracts $n/2$ from (\ref{eq:TotalSelectionBiasBCD}) to obtain the
excess selection bias in $n$ trials.
The average excess selection bias in $n$ trials (total excess selection
bias divided by $n$) of the $\mathit{BCD}(p)$ for different values of $n$ and
$p$ is provided in Table \ref{tab:SelectionBiasOfBCD}.

%
\begin{table}
\tablewidth=270pt
\caption{Average excess selection bias of the BCD for different values
of $n$ and $p$}
\label{tab:SelectionBiasOfBCD}
\begin{tabular*}{\tablewidth}{@{\extracolsep{\fill}}lcccc@{}}
\hline
& \multicolumn{4}{c@{}}{$\bolds p$} \\
& \multicolumn{4}{c@{}}{\hrulefill}\\
$\bolds n$ & \textbf{0.6} & \textbf{0.7}
& \textbf{0.8} & \textbf{0.9} \\
\hline
\phantom{00}5 & 0.058 & 0.107 & 0.146 & 0.177 \\
\phantom{0}10 & 0.070 & 0.129 & 0.178 & 0.217 \\
\phantom{0}15 & 0.072 & 0.129 & 0.173 & 0.207 \\
\phantom{0}20 & 0.075 & 0.136 & 0.183 & 0.220 \\
\phantom{0}25 & 0.076 & 0.135 & 0.179 & 0.213 \\
\phantom{0}50 & 0.080 & 0.140 & 0.186 & 0.221 \\
\phantom{0}75 & 0.081 & 0.140 & 0.185 & 0.219 \\
100 & 0.081 & 0.141 & 0.187 & 0.222 \\
200 & 0.082 & 0.142 & 0.187 & 0.222 \\
$\infty$ & 0.083 & 0.143 & 0.188 & 0.222 \\
\hline
\end{tabular*}
\end{table}

As expected, the excess selection bias increases with $p$. Also, note
that the average excess selection bias is not a monotonic function of
$n$. Asymptotically, the excess
is given in (\ref{selbias}) and is reported in the table under
$n=\infty$. One can see that the asymptotic formula is a good
approximation even for sample sizes as small as~50.

\subsection{Accidental bias}

With the help of Corollary \ref{crl:SigmaDn}, which provides the exact
form of the covariance matrix of the $\mathit{BCD}(p)$, one can compute the
accidental bias due to failure to adjust for any covariate $\mathbf{z}$,
given by $\mathbf{z}' \bolds{\Sigma}_{\mathbf{T}_n}
\mathbf{z}$. However, the point of the accidental bias is to
control the bias of the treatment effect caused by an \textit{unknown}
covariate. This leads to Efron's minimax solution of using the maximum
eigenvalue of $\bolds{\Sigma}_{\mathbf{T}_n}$, given in
inequality (\ref{minimax}). The maximum eigenvalue of $\bolds
{\Sigma}_{\mathbf{T}_n}$ therefore represents maximum
susceptibility to accidental bias. At this time we are able to prove
the following theorem.
\begin{theorem}
\label{thm:maxeigenvalbcd}
One of the eigenvalues of $\bolds{\Sigma}_{\mathbf
{T}_n}$ of the $\mathit{BCD}(p)$ is $2p$, for all $n \ge2$ and $p \in[1/2,1]$.
\end{theorem}
\begin{pf}
See Appendix \ref{appC}.
\end{pf}
\begin{Remark*}
The theorem affirms that the maximum eigenvalue of
$\bolds{\Sigma}_{\mathbf{T}_n}$ exceeds $1 + (p - q)^2$. This
shows that the maximum eigenvalue of the asymptotic covariance
structure studied by Efron (\citeyear{Efron1971}) and Steele
(\citeyear
{Steele1980}) is strictly less
than the maximum eigenvalue of $\bolds{\Sigma}_{\mathbf{T}_n}$.

We conjecture, based on vast numeric evidence, that the maximum
eigenvalue of $\bolds{\Sigma}_{\mathbf{T}_n}$ of the $\mathit{BCD}(p)$
does not depend on $n$ and is equal to $2p$ for all $n \ge2$ and $p
\in[1/2,1]$. Note that this leads to a maximum eigenvalue of 1 for
$p=0.5$, which is the maximum eigenvalue for complete randomization,
and 2 for $p=1$, which is the maximum eigenvalue for the permuted block
design with block size~2 [Rosenberger and Lachin (\citeyear
{RosenbergerLachin2002}),
Chapter 4].
\end{Remark*}

\subsection{Randomization tests}

The final application of these results is to random\-ization-based
inference procedures. Rosenberger and Lachin [(\citeyear
{RosenbergerLachin2002}), Chapters 7, 11]
discuss randomization tests in the context of linear rank statistics.
Let $\mathbf{Y}_n = (Y_1, Y_2,\ldots,Y_n)$ be the responses based
on some primary outcome variable, and let $\mathbf{y}_n$ be the
realization. The responses, $\mathbf{y}_n$, are treated as fixed
quantities, and under the randomization null hypothesis,
$\mathbf{y}_n$ is assumed to be unaffected by treatment assignments. The
observed difference between Groups $A$ and $B$ then only depends on the
manner the $n$ patients were randomized. The general form of linear
rank statistic is $W_n=\mathbf{a}_n'\mathbf{T}_n$ where
$\mathbf{a}_n = (a_{1n}, a_{2n},\ldots,a_{nn})'$ is a score
function of the ranks of $\mathbf{y}_n$. The scores $(a_{1n},
a_{2n},\ldots,a_{nn})'$ are usually centered by subtracting the
mean. Most standardly used test statistics in clinical trials have an
analogous formulation as a linear rank test.

Smythe and Wei (\citeyear{SW83}) and Hollander and Pe\~{n}a (\citeyear
{HP88}) noted that,
unlike for most other restricted randomization procedures, the test
$W_n$ is not asymptotically normal for the biased coin design.
Therefore the
computation of the test requires either the exact distribution or a
Monte Carlo approximation. While our results do not give the exact
distribution of the test statistic, we can compute its exact variance
as $\operatorname{Var}(W_n)=\mathbf{a}_n'\bolds{\Sigma}_{\mathbf
{T}_n}\mathbf{a}_n$ using Corollary \ref{crl:SigmaDn}. For example, using outcome
data from a diabetes trial given in Table 7.4 of Rosenberger and Lachin
(\citeyear{RosenbergerLachin2002}),
we generate a sequence of 50 treatment assignments from Efron's
$\mathit{BCD}(p=2/3)$ and obtain $W_n=-31$ with exact standard deviation
$100.52$. The latter computation required computing a $50 \times50$ matrix
using Corollary \ref{crl:SigmaDn}.

\section{Conclusions}\label{sec5}

Despite the favorable properties depicted in Efron's original paper,
the biased coin design is
sparsely used in clinical trials. The majority of clinical trials use a
permuted block design
which forces balance at regular intervals in the trial and achieves
perfect balance unless there
is an unfilled final block. However, in permuted blocks, some patients
are assigned to treatment with
probability 1 which can contribute to a vulnerability to selection
bias, particularly in unmasked trials.
We believe that Efron's procedure should be used regularly in clinical
trials where balance in
treatments is desirable, both for its simplicity and for the reason
that Efron suggested: it promotes balance with minimal susceptibility
to experimental biases. We now have quantified the distribution of
balance and the susceptibility to biases
in closed-form formulas for any $p$ and $n$, and this should aid the
clinical trialist in designing
the trial appropriately.

The selection of $p$ has always been an interesting question. Efron
used $p=2/3$ in some of his examples. At one extreme, $p=1/2$, we have
complete randomization which has minimal selection and accidental
biases, but maximum variability. At the other extreme, $p=1$, we have
a deterministic sequence with maximum selection and accidental biases,
but no variability. Formally, the selection
should be a trade-off between the degree of randomness desired (as
reflected in selection bias), accidental bias (which is linear) and
$\operatorname{Var}(D_n)$ which are competing objectives.
Such multi-objective problems can be solved through a compound
optimality criterion with weights reflecting the relative importance
of the criteria to the investigator. We now provide exact formulas for
these criteria in (\ref{eq:VarOfDn}) and (\ref{eq:TotalSelectionBiasBCD}).

We note that these results may have applicability beyond clinical
trials, as they form the basis of exact distribution theory
for a general asymmetric random walk. While the theorems are proved for
$p \ge0.5$, they can be generalized for any $p$ [Markaryan (\citeyear
{Markaryan09})].

\begin{appendix}
\section{\texorpdfstring{Proof of Theorem \protect\lowercase{\ref{thm:distdnshort}}}{Proof of Theorem 2.1}}\label{appA}

The following proposition follows immediately from the definition of
the $\mathit{BCD}(p)$ and is used without
explicit mention in the proof of Theorem \ref{thm:distdnshort}.
\begin{prp}
\label{prp:BasicPropertiesOfBCD}
Let $n = 1,2,3,\ldots,$ $k \in\mathbb{Z}$ and $q = 1-p$. The
following hold for the $\mathit{BCD}(p)$:
\begin{enumerate}
\item$P(D_{n} = k)>0 \iff|k| \leq n$ and $n$ and $k$ have the same parity;
\item$P(D_{n} = k) = P(D_{n} = -k)$;
\item$P(D_{n+1} = 0) = 2pP(D_{n} = 1)$;
\item$P(D_{n+1} = 1) = \frac{1}{2}P(D_{n} = 0) + pP(D_{n} = 2)$;
\item$P(D_{n+1} = k) = (1 - p)P(D_{n} = k - 1) + pP(D_{n} = k + 1)$,
for $2 \leq k \leq n$;
\item$P(D_{n+1} = n + 1) = (1 - p)P(D_{n} = n)$.
\end{enumerate}
\end{prp}

Next we formulate and prove two lemmas.
%
\begin{lem}
\label{lem:CaseK_1_l}
Let $n$ be a positive even integer, and let $l$ be an
integer satisfying $0 < l < n/2$. Then the following holds:
%
%
\begin{equation}
\label{eq:CaseK_1_l}\qquad
\frac{n-2l}{n+2l}
\pmatrix{
\dfrac{n}{2}+l \vspace*{2pt}\cr
l}
+
\frac{n-2l+4}{n+2l}
\pmatrix{
\dfrac{n}{2}+l \vspace*{2pt}\cr
l-1}
=
\frac{n+2-2l}{n+2+2l}
\pmatrix{
\dfrac{n}{2}+1+l \vspace*{2pt}\cr
l}.
\end{equation}
\end{lem}
\begin{pf}
First, we make a substitution, $u = n/2$ in (\ref{eq:CaseK_1_l}), to
obtain an equivalent expression,
%
%
\begin{equation}
\label{eq:CaseK_1_ls}
\frac{u-l}{u+l}
\pmatrix{
u + l \cr
l}+ \frac{u-l+2}{u+l}
\pmatrix{
u + l \cr
l - 1}
=
\frac{u + 1 - l}{u + 1 + l}
\pmatrix{
u + 1 + l \cr
l}.
\end{equation}
Using easily checked identities,
\[
\pmatrix{
u + l \cr
l}
=
\frac{u + 1}{u + 1 + l}
\pmatrix{u + 1 + l \cr l}
\]
and
\[
\pmatrix{u + l \cr
l - 1}
=
\frac{l}{u + 1 + l}
\pmatrix{
u + 1 + l \cr l},
\]
the left-hand side of (\ref{eq:CaseK_1_ls}) can be re-written as
\[
\frac{(u - l)(u + 1) + l(u - l + 2)}{(u + l)(u + 1 +l)}
\pmatrix{u + 1 + l \cr l}
\]
and the lemma follows from noting that
\[
\frac{(u - l)(u + 1) + l(u - l + 2)}{(u + l)} = u + 1 -l.
\]
\upqed\end{pf}
%
%
\begin{lem}
\label{lem:CaseK2_n}
Let $n$ be a positive integer, $k$ be an integer satisfying
$2 \leq k \leq n$, $l$ be an integer satisfying $1 \leq l \leq\frac{n
- k + 1}{2}$ and $n$ and $k$ have opposite parities. Then the following holds:
%
%
\begin{eqnarray}
\label{eq:CaseK2_n_l}
&&\frac{n + k -2l + 3}{n + k + 2l - 1}
\pmatrix{\dfrac{n + k + 1}{2} + l - 1 \cr
l - 1}
+ \frac{n + k - 2l - 1}{n + k + 2l - 1}
\pmatrix{\dfrac{n + k - 1}{2} + l \cr
l}
\hspace*{-32pt}\nonumber\\[-8pt]\\[-8pt]
&&\qquad =
\frac{n + k -2l + 1}{n + k + 2l + 1}
\pmatrix{\dfrac{n + k + 1}{2} + l \cr l}.
\hspace*{-32pt}\nonumber
\end{eqnarray}
\end{lem}
\begin{pf}
We first make a substitution, $u = (n + k + 1)/2$ in (\ref
{eq:CaseK2_n_l}), and obtain an equivalent expression,
%
%
\begin{equation}
\label{eq:CaseK2_n_ls}\qquad
\frac{u - l + 1}{u + l - 1}
\pmatrix{u + l - 1 \cr l - 1}
+
\frac{u - l - 1}{u + l - 1}
\pmatrix{u + l - 1 \cr l}
=
\frac{u - l}{u + l}
\pmatrix{u + l \cr l}.
\end{equation}
Using easily verified identities,
\[
\pmatrix{u + l - 1 \cr
l - 1}=
\frac{l}{u + l}
\pmatrix{u + l \cr l}
\quad\mbox{and}\quad
\pmatrix{u + l - 1\cr l}
=
\frac{u}{u + l}
\pmatrix{u + l \cr l}
\]
and dividing both sides of (\ref{eq:CaseK2_n_ls}) by
\[
\frac{1}{(u + l - 1)(u + l)}
\pmatrix{u + l \cr l},
\]
the result follows.
\end{pf}

Before we prove the theorem, note that in the light of Proposition \ref
{prp:BasicPropertiesOfBCD}, the assumptions on $n$ and $k$ are for the
purpose of identifying the nonzero probability events. Also, due to
symmetry, we can restrict the proof to the case of nonnegative~$k$. The
proof is by induction and involves a series of straightforward
calculations. The theorem is trivially true for the cases $n = 1$ and
$n = 2$.
We assume the theorem is true for all positive integers up to and
including $n$ and prove that it is true for $n + 1$.
The proof is broken out into four cases: $k = 0$, $k = 1$, $2 \leq k
\leq n$ and $k = n + 1$.
\begin{pf*}{Proof of Theorem \protect\ref{thm:distdnshort}}

\textit{Case $k = 0$.}
\begin{eqnarray*}
P(D_{n+1}= 0) &=& 2pP(D_{n} = 1) \\
&=&
2p \cdot
\frac{1}{2} p^{(n-1)/{2}}
\sum_{l=0}^{(n-1)/{2}}
{\frac{n+1-2l}{n+1+2l}
\pmatrix{\dfrac{n+1}{2}+l \cr
l}q^{1+l-1}
} \\
&=&
p^{(n+1)/{2}}
\sum_{l=0}^{(n-1)/{2}}
{\frac{n+1-2l}{n+1+2l}
\pmatrix{\dfrac{n+1}{2}+l \cr
l}q^{l}},
\end{eqnarray*}
which is exactly (\ref{eq:DistDnShortKzero}) with $n$ replaced by $n + 1$.

\textit{Case $k = 1$.} We need to show that
%
%
\begin{equation}
\label{eq:DistDnShortK1}
P(D_{n+1}= 1)= \frac{1}{2} p^{{n}/{2}}
\sum_{l=0}^{{n}/{2}}
{\frac{n+2-2l}{n+2+2l}
\pmatrix{\dfrac{n}{2}+1+l \cr
l}q^{l}}.
\end{equation}
Then
\begin{eqnarray*}
P(D_{n+1}= 1) &=& \frac{1}{2}P(D_{n} = 0) + pP(D_{n} = 2) \\
&=& \frac{1}{2}
p^{{n}/{2}}
\sum_{l=0}^{{n}/{2}-1}
{
\frac{n-2l}{n+2l}
\pmatrix{\dfrac{n}{2}+l \cr
l}q^{l}
}
\\
&&{} + p \cdot
\frac{1}{2} p^{({n-2})/{2}}
\sum_{l=0}^{({n-2})/{2}}
{\frac{n+2-2l}{n+2+2l}
\pmatrix{\dfrac{n+2}{2}+l \cr l}
q^{2+l-1}.}
\end{eqnarray*}
Now we shift the summation index in the second term, $l:=l + 1$, and
then collect the terms under a single summation,
%
%
\begin{eqnarray}
\label{eq:CaseK1_1}
&&P(D_{n+1}= 1) \nonumber\\
&&\qquad= \frac{1}{2}
p^{{n}/{2}}
\sum_{l=0}^{{n}/{2}-1}
{\frac{n - 2l}{n + 2l}
\pmatrix{\dfrac{n}{2} + l \cr
l}q^{l}
} \nonumber\\
&&\qquad\quad{}+
\frac{1}{2} p^{{n}/{2}}
\sum_{l=1}^{{n}/{2}}
{\frac{n + 2 - 2(l - 1)}{n + 2 + 2(l - 1)}
\pmatrix{\dfrac{n + 2}{2} + l - 1 \vspace*{2pt}\cr l - 1}
q^{l}} \\
&&\qquad=\frac{1}{2}
p^{{n}/{2}}
\Biggl\{\sum_{l=1}^{{n}/{2}-1}
{\left[\frac{n-2l}{n+2l}
\pmatrix{\dfrac{n}{2} + l \cr l}
+\frac{n-2l+4}{n+2l}
\pmatrix{\dfrac{n}{2} + l \vspace*{2pt}\cr
l - 1}
\right]}q^{l} \Biggr\} \nonumber\\
&&\qquad\quad{} + \frac{1}{2}
p^{{n}/{2}}
\left\{ 1 + \frac{2}{n}
\pmatrix{n \vspace*{2pt}\cr
\dfrac{n}{2} - 1}
q^{{n}/{2}}
\right\}. \nonumber
\end{eqnarray}
Similar to the right-hand side of (\ref{eq:DistDnShortK1}), the
expression obtained in (\ref{eq:CaseK1_1}) is a product of $p^{n/2}/2$
and a $(n/2)$th order polynomial in $q$. Therefore it remains to show
that the polynomial inside the curly braces in (\ref{eq:CaseK1_1}) is
the same as the polynomial in the right-hand side of (\ref{eq:DistDnShortK1}).
We will show term by term equality. First, the constant term in (\ref
{eq:CaseK1_1}) is $1$ which is the same as the constant term in (\ref
{eq:DistDnShortK1}). To show that the coefficients of $q^{{n}/{2}}$
are equal we need to show the following equality:
\[
\frac{2}{n}
\pmatrix{n \vspace*{2pt}\cr
\dfrac{n}{2} - 1}
=
\frac{1}{n + 1}
\pmatrix{n + 1\vspace*{2pt}\cr
\dfrac{n}{2}}.
\]
We transform the left-hand side to obtain the right-hand side as follows:
\begin{eqnarray*}
\frac{2}{n}
\pmatrix{n \vspace*{2pt}\cr
\dfrac{n}{2} - 1}
&=&
\frac{1}{{n}/{2}} \cdot\frac{n!}{({n}/{2}-1)!({n}/{2}+1)!}
=
\frac{n!}{({n}/{2})!({n}/{2} + 1)!} \\
&=&
\frac{1}{n + 1} \cdot\frac{(n + 1)!}{({n}/{2})!({n}/{2} + 1)!}
=
\frac{1}{n + 1}
\pmatrix{n + 1\vspace*{2pt}\cr
\dfrac{n}{2}}.
\end{eqnarray*}
To complete the proof for the case $k = 1$ it remains to show that the
coefficients of $q^l$ are equal for $0 < l < n/2$.
This is contained in Lemma \ref{lem:CaseK_1_l}.

\textit{Case $2 \leq k \leq n$.} We need to show that
%
%
\begin{eqnarray}
\label{eq:DistDnShortK2_n}\qquad
P(D_{n+1}= k) &=& \frac{1}{2} p^{({n - k + 1})/{2}}
\nonumber\\[-8pt]\\[-8pt]
&&{}\times
\sum_{l=0}^{({n - k + 1})/{2}}
{
\frac{n + k -2l + 1}{n + k + 2l + 1 }
\pmatrix{\dfrac{n + k + 1}{2} + l \cr
l}q^{k + l - 1}.}\nonumber
\end{eqnarray}
When $k = n$, $n + 1$ and $k$ have opposite parities; therefore we can
assume that $2 \leq k \leq n - 1$. We have
\begin{eqnarray*}
P(D_{n+1}= k) &=& pP(D_{n} = k + 1)
 + qP(D_{n} = k - 1) \\
&=& p \cdot\frac{1}{2}
p^{({n - k - 1})/{2}}
\sum_{l=0}^{({n - k - 1})/{2}}
{
\frac{n + k -2l + 1 }{n + k + 2l + 1}
\pmatrix{\dfrac{n + k + 1}{2} + l \cr l}
q^{k + l}} \\
&&{} + q \cdot\frac{1}{2}
p^{({n - k + 1})/{2}}\\
&&\hspace*{10pt}{}\times
\sum_{l=0}^{({n - k + 1})/{2}}
{
\frac{n + k -2l - 1}{n + k + 2l - 1}
\pmatrix{\dfrac{n + k - 1}{2} + l \cr
l}q^{k + l - 2}.}
\end{eqnarray*}
Now we shift the summation index in the first term, $l:=l + 1$, and
then collect the terms under a single summation to obtain
%
%
\begin{eqnarray}
\label{eq:DistDnShortK2_n_3}\hspace*{24pt}
&& P(D_{n+1}= k) \nonumber\\
&&\qquad=\frac{1}{2} p^{({n - k + 1})/{2}} q^{k - 1}
\sum_{l = 1}^{({n - k + 1})/{2}}
\frac{n + k -2l + 3}{n + k + 2l - 1}
\pmatrix{\dfrac{n + k + 1}{2} + l - 1 \vspace*{2pt}\cr
l - 1}
q^{l} \nonumber\\
&&\qquad\quad{} + \frac{1}{2}
p^{({n - k + 1})/{2}} q^{k - 1}
\sum_{l = 0}^{({n - k + 1})/{2}}
\frac{n + k - 2l - 1}{n + k + 2l - 1}
\pmatrix{\dfrac{n + k - 1}{2} + l \vspace*{2pt}\cr
l}q^{l}
\\
&&\qquad = c \Biggl\{
\sum_{l = 1}^{({n - k + 1})/{2}}
\Biggl[
\frac{n + k -2l + 3}{n + k + 2l - 1}
\pmatrix{\dfrac{n + k + 1}{2} + l - 1 \vspace*{2pt}\cr
l - 1} \Biggr] \Biggr\}
\nonumber\\
&&\qquad\quad{} +
c \Biggl\{
\sum_{l = 1}^{({n - k + 1})/{2}}
\Biggl[\frac{n + k - 2l - 1}{n + k + 2l - 1}
\pmatrix{\dfrac{n + k - 1}{2} + l \vspace*{2pt}\cr
l} \Biggr] q^{l} +1 \Biggr\}, \nonumber
\end{eqnarray}
where $c = p^{(n-k+1)/2} q^{k - 1}/2$. Comparing (\ref
{eq:DistDnShortK2_n_3}) with (\ref{eq:DistDnShortK2_n}) we immediately
see that the terms corresponding to $l = 0$ are equal to $c$. To
complete the proof for the case $2 \leq k \leq n$ all that remains is
an application of Lemma \ref{lem:CaseK2_n}.

\textit{Case} $k = n + 1$: This follows immediately from the fact that
\[
P(D_{n+1}= n+1)= \tfrac{1}{2} q^n.
\]

The theorem is proved.
\end{pf*}\vspace*{-15pt}

\section{\texorpdfstring{Proof of Theorem \protect\lowercase{\ref{thm:distdndm}}}{Proof of Theorem 2.2}}\label{appB}

The following proposition follows immediately from the Markovian
property and time homogeneity of the $\mathit{BCD}(p)$ process.
\begin{prp}
\label{prp:BasicMarkovianPropertiesOfBCD}
Let $n = 0,1,2,3,\ldots,$ $m = 1,2,3,\ldots$ and $m \ge n $.
Define $\sigma(\mathbf{T}_n)$ to be the sigma-algebra generated by
$T_1,\ldots,T_n$. The following hold for the $\mathit{BCD}(p)$:
\begin{enumerate}
\item$P(T_{m} = \pm1 | D_{n}, \sigma(\mathbf{T}_{n})) = P(T_{m} =
\pm1 | D_{n})$;
\item$P(T_{m} = \pm1 | D_{n} = k) = P(T_{m + l} = \pm1 | D_{n + l} =
k)$, for any $l \ge0$.
\end{enumerate}
\end{prp}

Next we state and prove three lemmas that are used in the
proof of Theorem \ref{thm:distdndm}.
%
%
\begin{lem}
\label{lem:Joint_To_Conditional}
Let $1 \leq n < m$. Then the following holds for the
$\mathit{BCD}(p)$, $p \in[0,1]$:
%
%
\begin{equation}
\label{eq:Joint_To_Conditional}
P(T_{n}=1, T_{m}=1) =
\sum_{k=-n+1}^{n-1}
{
P(T_{m} = 1 | D_{n} = k + 1)
d_{n-1,k} t_{k}.
}
\end{equation}
\end{lem}
\begin{pf}
Before providing the proof, note that because Theorem \ref
{thm:distdnshort} gives the form of $d_{n,k}$, the lemma reduces the
finding of $P(T_{n}=1, T_{m}=1)$ to finding conditional probabilities
of the form $P(T_{m} = 1 | D_{n} = k)$.
By conditioning on $D_{n-1}$ we obtain
%
%
\begin{eqnarray}
\label{eq:TmTn_Cond_on_Dn}
&&P(T_{n}=1, T_{m}=1) =
\sum_{k=-n+1}^{n-1}
P(T_{n}=1, T_{m}=1 | D_{n-1} = k)\nonumber\\[-8pt]\\[-8pt]
&&\hspace*{128.1pt}{} \times P(D_{n- 1} = k).\nonumber
\end{eqnarray}
Note that (\ref{eq:TmTn_Cond_on_Dn}) holds for $n = 1$ as
well because we defined $P(D_{0} = 0) = 1$. Also, instead of requiring
$n - k$ be odd so that $P(D_{n-1} = k) > 0$, we follow the adopted
convention that the probabilities of events conditional on
zero-probability event are treated as $0$.

Now we make use of an easily verified identity,\setcounter{footnote}{1}
\[
P(A \cap B | C) = P(A | B \cap C) \cdot P(B | C)\mbox{,}\footnote{This
identity still holds when either $B$ or $C$ have zero probability when
used with the adopted convention.}
\]
to transform the conditional probabilities in the right-hand side of
(\ref
{eq:TmTn_Cond_on_Dn}),
%
\begin{eqnarray}
\label{eq:Eval_Cond_Prob}
&&P(T_{n}=1, T_{m}=1 | D_{n-1} = k)
\nonumber\\[-8pt]\\[-8pt]
&&\qquad
= P(T_{m}=1 | T_{n}=1, D_{n-1} = k)P(T_n = 1| D_{n-1} = k).\nonumber
\end{eqnarray}
Now we use the fact that the following two events are equal:
\[
\{D_{n-1} = k \mbox{ and } T_n = 1\} \quad\mbox{and}\quad \{D_{n} = k+1 \mbox{
and } T_n = 1\}
\]
and that $T_m$ is conditionally independent of $T_n$ given $D_n$ to write
\begin{eqnarray*}
P(T_{m}=1 | T_{n}=1, D_{n-1} = k) &=& P(T_{m}=1 | T_{n}=1, D_{n} =
k+1) \\
&=& P(T_{m}=1 | D_{n} = k+1).
\end{eqnarray*}
Substituting this last expression into the right-hand side of (\ref
{eq:Eval_Cond_Prob}) we obtain
\[
P(T_{n}=1, T_{m}=1 | D_{n-1} = k) = P(T_{m}=1 | D_{n} = k+1)P(T_n = 1|
D_{n-1} = k).
\]
The result follows from substitution into (\ref{eq:TmTn_Cond_on_Dn}).
\end{pf}

The next lemma is devoted to finding the first visit probabilities of
the imbalance process into the 0 state. We define
$\tau_i$ to be the number of steps the imbalance process makes to visit
state $0$ for the first time from the
$i$th state.
%
%
\begin{lem}
\label{lem:First_Visit_Probabilities}
For the imbalance process of the $\mathit{BCD}(p)$, $p \in[0,1]$,
the probabilities of the first visits from state $k$, $k = \pm1, \pm
2,\ldots,$ into state $0$ in exactly $l$ steps, $l \ge|k|$, is given by
the following formula:
%
%
\begin{equation}
\label{eq:FirstVisitProbabilities}
f_{k,0}^{(l)} = P(\tau_{k} = l) =
\frac{|k|}{l}
\pmatrix{l \vspace*{2pt}\cr
\dfrac{l+|k|}{2}}
p^{({l+|k|})/{2}} q^{({l-|k|})/{2}},
\end{equation}
where, according to the adopted convention, the combination is to be
treated as $0$ when $(l+|k|)/2$ is not an integer.
\end{lem}
\begin{pf}
First, due the symmetry, $ f_{k,0}^{(l)} = f_{-k,0}^{(l)}$, for any $k
\in\mathbb{N}$. Therefore without loss of generality, we can assume
that $k$ is positive. Thus we are concerned with finding first visit
probabilities from state $k, k \in\mathbb{N}$ into state $0$ in
exactly $l$, $l \ge k$, steps.

We can treat this problem as a random walk on the nonnegative integers
with an absorbing barrier at $0$ and use well-known results in the
classical gambler's ruin problem where the gambler plays with
infinitely reach adversary and at each step wins one unit with
probability $q$ and loses one unit with probability $p$. The question
is equivalently formulated as: what is the probability that a gambler
with initial capital of $k$, $k \in\mathbb{N}$, is ruined in exactly
$l$, $l \ge k$, steps? These probabilities are well known and can be
found in (4.14) of Feller (\citeyear{Feller1968}). One needs
to reverse the
roles of $p$ and $q$ and replace $z$ with $k$ and $n$ with $l$.
\end{pf}

Lemma \ref{lem:First_Visit_Probabilities} provides all the
nontrivial probabilities for $f_{k,0}^{(l)}$. To complete the remaining
cases, we note that $f_{k,0}^{(0)} = 0$ when $k \ne0$, and
$f_{0,0}^{(0)} = 1$.

The next lemma provides probabilities for the imbalance process to
ultimately reach the 0 state from any other state.
%
%
\begin{lem}
\label{lem:Evantual_Visiting_Probabilities}
For the imbalance process of the $\mathit{BCD}(p)$, $p \in[0.5,1]$,
the probability of ultimately reaching state $0$ from state k, $k = \pm
1, \pm2,\ldots,$ is 1.
\end{lem}
\begin{pf}
The proof of the lemma is similar to that of Lemma \ref
{lem:First_Visit_Probabilities}. Again, without loss of generality, it
can be assumed that $k$ is positive. The problem is equivalent to
computing the probability of ultimate ruin in the classical gambler's
ruin problem when the gambler, having an initial capital $k$, plays
with infinitely reach adversary and at each step wins one unit with
probability $q$ and loses one unit with probability $p$. These
probabilities can be found in (2.18) of Feller (\citeyear
{Feller1968}). One
needs to reverse the roles of $p$ and $q$ and replace $z$ with $k$.
\end{pf}

Note that Lemma \ref{lem:Evantual_Visiting_Probabilities} implies that $f_{k,0}^{(l)}$ is a
probability mass function when $p \in[0.5,1]$.

Before starting the proof of Theorem \ref{thm:distdndm}, note that only about half of
the summands in the right-hand side of (\ref{eq:JointTnTm}) will be nonzero
because $d_{n,k} = 0$ whenever $n-k$ is not even.
\begin{pf*}{Proof of Theorem \protect\ref{thm:distdndm}}
The essence of the proof is in evaluating conditional probabilities of
the form $P(T_m = 1 | D_n = k)$. We will show that for $1 \le n < m$,
$|k| \le n$ and $n - k$ even, the following holds:
%
%
\begin{equation}
\label{eq:P_Tm_Given_Dn_k}
P(T_m = 1 | D_n = k) =
\bigl(\tfrac{1}{2}-t_{k}
\bigr)
\hat{f}_{k,0}^{(m-n-1)} + t_{k}.
\end{equation}
This equation is of independent interest as it provides the form of
probabilities of treatment assignments conditional on a past value of imbalance.
Note that when $k = 0$, (\ref{eq:P_Tm_Given_Dn_k}) simply
states that $P(T_m = 1 | D_n = 0) = 1/2$, as expected. The case when $m
= n + 1$ is the definition of the $\mathit{BCD}(p)$.

To prove (\ref{eq:P_Tm_Given_Dn_k}), we use a conditioning argument and
condition on the first visit events into the 0 state, as follows:
%
%
\begin{eqnarray}
\label{eq:Cond_on_Tau}
&&P(T_{m}= 1 | D_n = k)\nonumber\\
&&\qquad=\sum_{l = 0}^{m - n - 1}P(T_m = 1 | D_n = k,
\tau
_{k} = l)P(\tau_{k} = l | D_n = k)
\nonumber\\[-8pt]\\[-8pt]
&&\qquad\quad{}
+P(T_m = 1 | D_n = k, \tau_{k} \notin[0, m-n-1])\nonumber\\
&&\qquad\quad\hspace*{10.3pt}{}\times P(\tau_{k} \notin[0,
m-n-1] | D_n = k).\nonumber
\end{eqnarray}
We first evaluate $P(T_m = 1 | D_n = k, \tau_{k} = l)$ for
the case $0 \le l \le m - n - 1$.
We only need to look at the cases when $(l - |k|)/2$ is a nonnegative
integer because in all other cases $P(\tau_{k} = l | D_n = k) = 0$:
\begin{eqnarray*}
&& P(T_m = 1 | D_n = k, \tau_{k} = l)\\
&&\qquad =P(T_m = 1 | D_n = k, D_{n + 1} \ne0, D_{n + 2} \ne0,\ldots,
D_{n +
l - 1} \ne0, D_{n + l} = 0) \\
&&\qquad =P(T_m = 1 | D_{n + l} = 0) =P(T_{m - n - l} = 1) =1/2.
\end{eqnarray*}
The first equality, in the chain of equalities above, is a
consequence of the following equality of events:
\[
\{D_n = k, \tau_{k} = l\} = \{D_n = k, D_{n + 1} \ne0, D_{n + 2} \ne
0,\ldots, D_{n + l - 1} \ne0, D_{n + l} = 0\}.
\]
The second equality is just the Markovian property of the
imbalance process [see Proposition~\ref
{prp:BasicMarkovianPropertiesOfBCD}(1)].
The third equality follows from time-homogeneity property formulated in
Proposition \ref{prp:BasicMarkovianPropertiesOfBCD}(2).
Thus we have proved that when $1 \le n < m$, $|k| \le n$, $n - k$ is
even, $0 \le l \le m - n - 1$ and $(l - |k|)/2$ is a nonnegative
integer, then
%
%
\begin{equation}
\label{eq:Case_l_smaller_m-n-1}
P(T_m = 1 | D_n = k, \tau_{k} = l) = 1/2.
\end{equation}
Now we turn to the case when $\tau_{k} \notin[0, m-n-1]$. As
before, we have $1 \le n < m$, $|k| \le n$ and $n - k$ is even. We look
at three sub-cases.

\textit{Case $k>0$.}
\begin{eqnarray*}
&&P(T_m = 1 | D_n = k, \tau_{k} \ge m - n)\\
&&\qquad =P(T_m = 1 | D_n = k, D_{m - 1} > 0, \tau_{k} \notin[0, m-n-1]) \\
&&\qquad =P(T_m = 1 | D_{m - 1} > 0)=q.
\end{eqnarray*}
The first equality above follows from equality of
$\{D_n = k, \tau_{k} \ge m - n)\}$ and
$\{D_n = k, D_{m - 1} > 0, \tau_{k} \notin[0, m-n-1]\}$.
The second equality follows from Proposition~\ref
{prp:BasicMarkovianPropertiesOfBCD}(1).

\textit{Case $k < 0$.}
\begin{eqnarray*}
&& P(T_m = 1 | D_n = k, \tau_{k} \ge m - n)\\
&&\qquad =P(T_m = 1 | D_n = k, D_{m - 1} < 0, \tau_{k} \notin[0, m-n-1]) \\
&&\qquad =P(T_m = 1 | D_{m - 1} < 0)=p.
\end{eqnarray*}
The first equality above follows from equality of
$\{D_n = k, \tau_{k} \ge m - n)\}$ and
$\{D_n = k, D_{m - 1} < 0, \tau_{k} \notin[0, m-n-1]\}$.
The second equality follows from Proposition~\ref
{prp:BasicMarkovianPropertiesOfBCD}(1).

\textit{Case $k = 0$.}
\[
P(T_m = 1 | D_n = 0, \tau_{0} \notin[0, m-n-1]) = 0,
\]
because of impossibility of the event $\{\tau_{0} \ge1\}$.

Substituting (\ref{eq:Case_l_smaller_m-n-1}) and the expressions
obtained in the above three cases into (\ref{eq:Cond_on_Tau}), we obtain
%
%
\begin{equation}
\label{eq:Tm_Cond_Dn_Part1}\qquad
P(T_m = 1 | D_n = k) = \sum_{l = 0}^{m - n - 1}
{\frac{1}{2} f_{k,0}^{(l)} + t_{k} P(\tau_{k} \notin[0, m-n-1] | D_n
= k).}
\end{equation}
According to Lemma \ref{lem:Evantual_Visiting_Probabilities}, when $p
\in[1/2, 1]$, we have
%
%
\begin{equation}
\label{eq:Tm_Cond_Dn_Part2}\qquad
P(\tau_{k} \notin[0, m-n-1] | D_n = k) = 1 - P(\tau_{k} \in[0,
m-n-1] | D_n = k).
\end{equation}
Substituting (\ref{eq:Tm_Cond_Dn_Part2}) into (\ref
{eq:Tm_Cond_Dn_Part1}), we obtain
\begin{eqnarray*}
P(T_m = 1 | D_n = k) & = & \sum_{l = 0}^{m - n - 1}
{\frac{1}{2} f_{k,0}^{(l)} + t_{k} \bigl(1 - \hat{f}_{k,0}^{(m-n-1)} \bigr)
} \\
& = & \frac{1}{2} \hat{f}_{k,0}^{(m - n -1)} + t_{k}
\bigl(1 - \hat{f}_{k,0}^{(m-n-1)} \bigr) \\
& = & \biggl(\frac{1}{2} - t_{k} \biggr) \hat{f}_{k,0}^{(m -
n -1)} + t_{k}.
\end{eqnarray*}
Thus (\ref{eq:P_Tm_Given_Dn_k}) is proved. To complete the
proof of the theorem, it remains to use Lemma \ref
{lem:Joint_To_Conditional} and substitute (\ref{eq:P_Tm_Given_Dn_k})
with $k:=k+1$ into (\ref{eq:Joint_To_Conditional}).
\end{pf*}

\section{\texorpdfstring{Proof of Theorem \protect\lowercase{\ref{thm:maxeigenvalbcd}}}{Proof of Theorem 4.1}}\label{appC}

We will show that $2p$ is an eigenvalue of $\bolds{\Sigma
}_{\mathbf{T}_n}$ with the corresponding normalized eigenvector
$\mathbf{a}^n = (a_1^n,a_2^n,a_3^n, \ldots, a_n^n)' = (\sqrt{2}/2,
-\sqrt{2}/2,0, \ldots, 0)'$. The proof proceeds by induction. The
theorem is trivially true for the case $n = 2$. For the case $n = 2$,
one can actually show that $2p$ is the maximum eigenvalue.\footnote{The
same can be shown for $n = 3$ and $n = 4$ by solving for the zeroes of
the characteristic polynomials.} The two eigenvalues of $\bolds
{\Sigma}_{\mathbf{T}_2}$ are $2p$ and $2 - 2p$ and $2p \ge2 - 2p$
when $p \ge1/2$.
\begin{pf}
We assume the theorem is true for all positive integers $n \ge2$, and
prove that it is true for $n + 1$. We partition $\bolds{\Sigma
}_{\mathbf{T}_{n + 1}}$ as follows:
\[
\left[
\begin{array}{c | c}
\bolds{\Sigma}_{\mathbf{T}_n} & \mathbf{b} \\ \hline
\mathbf{b}' & 1
\end{array}
\right],
\]
where $\mathbf{b} = (\sigma_{1,n+1}, \sigma_{2,n+1}, \ldots, \sigma
_{n,n+1})'$ with $\sigma_{ij} = \operatorname{Cov}(T_i, T_j)$. Denote
\[
\mathbf{x} = (\mathbf{x}_1, 0)',
\]
where $\mathbf{x}_1$ is the $n$-dimensional vector,
\[
\mathbf{x}_1 = \bigl(\sqrt{2}/2, -\sqrt{2}/2,0, \ldots, 0\bigr)'.
\]
We need to show that
%
%
\begin{equation}
\label{eq:Eigen_BCD_n_plus_1}
\left[
\begin{array}{c | c}
\bolds{\Sigma}_{\mathbf{T}_n} & \mathbf{b} \\ \hline
\mathbf{b}' & 1
\end{array}
\right]
\left[
\begin{array}{c}
\mathbf{x}_1 \\ \hline
0
\end{array}
\right]
=
2p
\left[
\begin{array}{c}
\mathbf{x}_1 \\ \hline
0
\end{array}
\right].
\end{equation}
By the induction assumption, we have that $\bolds{\Sigma
}_{\mathbf{T}_n} \mathbf{x}_1 = 2p \mathbf{x}_1$. To prove
(\ref{eq:Eigen_BCD_n_plus_1}), it remains to show that
$\mathbf{b}' \mathbf{x}_1 = 0$. This is equivalent to $\sqrt
{2}/2(\sigma_{1,n+1} - \sigma_{2,n+1}) = 0$ which in turn is equivalent
to (see Corollary \ref{crl:SigmaDn})
%
%
\begin{equation}
\label{eq:T1Tn1=T2Tn1}
P(T_1 = 1, T_{n + 1} = 1) = P(T_2 = 1, T_{n + 1} = 1).
\end{equation}
From Theorem \ref{thm:distdndm} we have the forms of $P(T_1 =
1, T_{n + 1} = 1)$ and $P(T_2 = 1, T_{n + 1} = 1)$,
\begin{eqnarray*}
P(T_{1}=1, T_{n + 1}=1) & = &
\bigl(
\bigl(
\tfrac{1}{2} - q
\bigr)
\hat{f}_{1,0}^{(n - 1)} + q
\bigr)
\tfrac{1}{2}
=
\tfrac{1}{2}
\bigl(
\tfrac{1}{2} - q
\bigr)
\hat{f}_{1,0}^{(n - 1)}
+
\tfrac{1}{2}q, \\
P(T_{2}=1, T_{n + 1}=1) & = &
\bigl(
\bigl(
\tfrac{1}{2} - \tfrac{1}{2}
\bigr)
\hat{f}_{0,0}^{(n - 2)} + \tfrac{1}{2}
\bigr)
\tfrac{1}{2}p
+
\bigl(
\bigl(
\tfrac{1}{2} - q
\bigr)
\hat{f}_{2,0}^{(n - 2)} + q
\bigr)
\tfrac{1}{2}q \\
& = &
\tfrac{1}{4}p + \tfrac{1}{2}q
\bigl(
\tfrac{1}{2} - q
\bigr)
\hat{f}_{2,0}^{(n - 2)}
+
\tfrac{1}{2}q^2.
\end{eqnarray*}
Thus, in order to show (\ref{eq:T1Tn1=T2Tn1}), we need to
show that
%
%
\begin{equation}
\label{eq:T1Tn1=T2Tn1_detail}
\tfrac{1}{2}
\bigl(
\tfrac{1}{2} - q
\bigr)
\hat{f}_{1,0}^{(n - 1)}
+
\tfrac{1}{2}q
=
\tfrac{1}{4}p + \tfrac{1}{2}q
\bigl(
\tfrac{1}{2} - q
\bigr)
\hat{f}_{2,0}^{(n - 2)}
+
\tfrac{1}{2}q^2.
\end{equation}
Using an easily verified identity,
\[
\hat{f}_{1,0}^{(n - 1)} = p + q \hat{f}_{2,0}^{(n - 2)}
\]
and substituting it into (\ref{eq:T1Tn1=T2Tn1_detail}), we obtain
%
%
\begin{eqnarray}
\label{eq:T1Tn1=T2Tn1_detail2}
&&\tfrac{1}{2}
\bigl(
\tfrac{1}{2} - q
\bigr)
\bigl(
p + q \hat{f}_{2,0}^{(n - 2)}
\bigr)
+
\tfrac{1}{2}q
\nonumber\\[-8pt]\\[-8pt]
&&\qquad=
\tfrac{1}{4}p + \tfrac{1}{2}q
\bigl(
\tfrac{1}{2} - q
\bigr)
\hat{f}_{2,0}^{(n - 2)}
+
\tfrac{1}{2}q^2.\nonumber
\end{eqnarray}
The term $\frac{1}{2}q (\frac{1}{2} - q ) \hat
{f}_{2,0}^{(n - 2)}$ appears in both sides of (\ref
{eq:T1Tn1=T2Tn1_detail2}); subtracting it from both sides we require
\[
\tfrac{1}{2}p
\bigl(
\tfrac{1}{2} - q
\bigr)
+
\tfrac{1}{2}q
=
\tfrac{1}{4}p + \tfrac{1}{2}q^2.
\]
This last equality is trivially checked.
\end{pf}
\end{appendix}

\printaddresses

\end{document}